\advance\year by -2000

\def\os{\obeyspaces}
\overfullrule=0pt\baselineskip=24pt
\def\C{{\rm C\hskip -6pt \vrule \hskip
6pt}}

\def\R{\os{\bf{R}\os}}

\def\s{\os{\bf{s}\os}}

\def\C{\os{\bf{C}\os}}
\def\D{\os{\bf{D}\os}}

\def\v{\os{\bf{v}\os}}
\def\r{\os{\bf{r}\os}}
\def\w{\os{\bf{w}\os}}
\def\a{\os{\bf{a}\os}}
\def\b{\os{\bf{b}\os}}
\def\s{\os{\bf{s}\os}}
\def\S{\os{\bf{S}\os}}
\def\G{\os{\bf{G
}\os}}
\def\x{\os{\bf{x}\os}}
\def\G{{\bf G\obeyspaces}}

\def\H{\os{\bf{H}\os}}
\def\T{\os{\bf{T}\os}}
\def\CC{\os{\bf{CC}\os}}
\def\CD{\os{\bf{CD}\os}}
\def\DC{\obeyspaces{\bf{DC}\os}}
\def\DD{\obeyspaces{\bf{DD}\os}}
\def\HH{\os{\hbox{$\H\otimes\H$}}\os}

\def\TT{\os{\hbox{$\T\otimes\T$}}\os}

\def\R{{\bf R }}
\def\qed{\hfill\vrule width 5pt height 5pt}

\font\titlefont=cmr10 at 17pt
\def\p{\os{\bf{p}\os}}
\def\q{\os{\bf{q}\os}}

\def\u{\os{\bf{u}\os}}
\def\v{\os{\bf{v}\os}}
\def\uv{\os{\bf{uv}\os}}
\def\vu{\os{\bf{vu}\os}}

\def\pu{\os{\bf{pu}\os}}
\def\pv{\os{\bf{pv}\os}}
\def\pvu{\os{\bf{pvu}\os}}

\line{\hfill
\the\month/\the\day/0\the\year}

\centerline{\titlefont Nash Equilibria
in Quantum Games}
\centerline{by}
\centerline{Steven E. Landsburg}
\centerline{University of Rochester}
\bigskip

\baselineskip=12pt

\bigskip

\baselineskip=24pt
{\bf Introduction.}

Quantum game theory models the behavior of 
strategic agents with access to quantum technology.
For example, agents might use observations of 
entangled particles to randomize their strategies
as in [CHTW] or [DL]; alternatively they might 
use quantum devices to communicate with each other
or with the referee.

Among models of quantum communication, the most
widely studied is the {\it EWL model\/} ([EW], [EWL])
of Eisert, Wilkens and Lewenstein.  Here we start
with an ordinary game \G, envision players who
communicate with a referee via a specific quantum
protocol, and, motivated by this vision,
 construct a new game $\G^Q$ with greatly enlarged
strategy spaces.

If \G\phantom{,} is a two-by-two game (that is:  a game with 
two players, each of whom has a two point strategy
space), then the game 
$\G^Q$ has strategy spaces that are naturally
identified with the 3-sphere $\S^3$.  A mixed strategy
in the game $\G^Q$ is therefore an arbitrary
 probability distribution on
$\S^3$.  

Partly because the space of mixed strategies is so
large, mixed-strategy Nash equilibria in the game 
${\bf G}^Q$ can be difficult to find.  The 
present paper ameliorates this difficulty by proving
that up to a natural notion of equivalence, all Nash equilbria
have a particularly simple form.  More specifically (and still
up to equivalence, as defined in Section 1):

1)  Every mixed strategy that occurs in a Nash equilibrium 
is supported on at most four points of ${\bf S}^3$.

2)  Those four points must lie in one of a small number of
allowable geometric configurations.

I call a two-by-two game ``generic'' if, for each player, the four 
possible payoffs are all distinct and the six pairwise sums of those 
payoffs are all distinct. In this paper, I will state 
and prove the main theorem for generic games, referring 
the reader to my unpublished working paper [NE] for the 
(considerably uglier but no more enlightening) generalization 
to the non-generic case.

Section 1 lays out the motivation and the details of the EWL 
model.  Section 2 presents the main technical lemmas.  
Section 3 contains the main classification theorem 
(3.3). 
Section 4 addresses some natural questions raised by the 
statement of the main theorem.   Section 5 collects 
a few additional remarks and applications; 
the most striking is that in any mixed strategy 
quantum equilibrium of any two-by-two zero sum game, 
each 
player earns exactly the average of the four possible 
payoffs.

{\bf 1.  The Eisert-Wilkens-Lewenstein Model.}

Let \G be a two-player game with strategy sets ${\cal S}_1$, 
${\cal S}_2$.  The EWL construction is motivated by specifying a 
purely classical communication
protocol:  A referee issues each player a penny in in one of
two states \H (``heads up'') or \T (``tails up'').  A player
indicates his choice of strategy by returning his penny  either
flipped or unflipped.  The referee observes the returned pennies
and computes payoffs accordingly.

Now replace the pennies with subatomic particles whose state
spaces are complex vector spaces with basis $\{\H,\T\}$.  A state
is a line through the origin, which we will often denote by 
specifying some nonzero point on that line.  The state space for the pair of
particles is the tensor product of the two individual state spaces.

The referee prepares two pennies in the {\it maximally entangled
state\/} $\H\otimes\H+\T\otimes\T$.  Each player returns the penny
after acting on its state by the special unitary operator of his
choice.  The two classical strategies (i.e. the strategies in the 
game \G) are identified with the operators
$$\C=\pmatrix{1&0\cr0&1}\qquad\D=\pmatrix{0&1\cr -1&0\cr}\eqno(1.0.1)$$

If Players One and Two select the unitary operators {\bf U} and {\bf V},
we denote the resulting state
by
$${\bf UV}=({\bf U}\otimes 1)\phantom{n}(\HH+\TT)
\phantom{n}(1\otimes{\bf V})\eqno(1.0.2)$$
The referee then performs an observation with the four possible outcomes
\CC, \CD, \DC and \DD, and makes payoffs accordingly; if we write
$${\bf UV}=
\alpha_1\CC+\alpha_2\DD+\alpha_3\CD+\alpha_4\DC\eqno(1
.0.3)
$$
(with complex scalar coefficients),  then the 
probabilities 
of the four states are proportional to 
$|\alpha_1|^2,|\alpha_2|^2,|\alpha_3|^2,|\alpha_4|^2$.

Now identify Player One's strategy space with the unit 
quaternions by mapping the unitary matrix with top row $(A,B)$
to the quaternion $A+Bj$; identify Player Two's
strategy space with the unit quaternions by mapping
the unitary matrix with top row $(P,Q)$ to the quaternion $P-
jQ$.
From (1.0.2) one readily calculates the coefficients in 
(1.0.3) 
and discovers 
the following remarkably simple formula:

{\bf Proposition 1.1.}  {\sl Suppose Player One plays
the quaternion \p and Player Two plays the quaternion
\q.  Then for $t=1,\ldots,4$, we have 
$$|\alpha_t|=2\Big|\pi_t(\p\q )\Big|$$
where the $\pi_t$ are the coordinate functions defined 
by 
$$\p=\pi_1(\p)+\pi_2(\p)i+\pi_3(\p)j+\pi_4(\p)k$$}

Motivated by Proposition 1.1 and the preceding 
discussion, we 
make the following definitions:

{\bf Definitions and Remarks 1.2.}  Let \G be a two by 
two 
game with 
strategy spaces $S_i=\{\C,\D\}$ and payoff functions 
$P_i:S_1\times S_2\rightarrow \R$.  Then the {\sl 
associated 
quantum game\/} $\G^Q$ is the two-player game in which 
each 
strategy space is the unit quaternions, and payoffs are 
calculated as 
$$P_i^Q(\p,\q)=\pi_1(\p\q)^2P_i(\C,\C)+\pi_2(\p\q)^2P_i(
\D,\D
)+
\pi_3(\p\q)^2P_i(\C,\D)+\pi_4(\p\q)^2P_i(\D,\C)$$

Note that for any strategy $\p$ chosen by Player 1, and 
for 
any probability distribution whatsoever over the four 
outcomes $(\C,\C)$, etc., Player 2 can always adopt a 
strategy $\q$ that effects this probability 
distribution: Let 
$a^2,b^2,c^2,d^2$ be the desired probabilities, let 
$\r=a+bi+cj+dk$ and set $\q=\p^{-1}\r$.  Therefore, in 
the 
game $\G^Q$, there can never be an equilibrium in pure 
strategies unless one of the four strategy pairs leads 
to an 
optimal outcome for both players.
Thus in $\G^Q$, pure-strategy equilibria are both rare 
and 
uninteresting.   

Next we consider mixed strategies.   
A {\sl mixed quantum strategy\/} for \G is a mixed 
strategy 
in the game $\G^Q$, i.e. a probability distribution on 
the 
space of unit quaternions.  If \p is a unit quaternion, 
I will sometimes identify 
\p with the mixed strategy supported entirely on \p.  If 
$\nu$ 
and $\mu$ are mixed strategies, I will write 
$P_i^Q(\nu,\mu)$ 
for the corresponding expected playoff to player $i$; 
that 
is:
$$P_i^Q(\nu,\mu)=\int P_i(\p,\q)d\nu(\p)d\mu(\q)$$

This gives rise to a new game $\G^Q_{\hbox{\bf mixed}}$, in which the
strategy sets are the sets of mixed quantum strategies in $\G$.  
A {\it mixed-strategy
Nash equilibrium\/} in $\G^Q$ is a Nash equilibirium in 
$\G^Q_{\hbox{\bf mixed}}.$

Our goal is to classify the mixed-strategy Nash equilibria in $\G^Q$.
The 
definitions that occupy the remainder of this section 
kick 
off that process by partitioning the set of Nash 
equilibria 
into natural equivalence classes.

{\bf Definition 1.3.}  Two mixed strategies $\mu$ and 
$\mu'$ 
are {\sl equivalent\/} if
$$\int\pi_t(\p\q)d\mu(\q)= \int\pi_t(\p\q)d\mu'(\q)$$
for all unit quaternions $\p$ and all $t=1,2,3,4$.

In other words, $\mu$ and $\mu'$ are equivalent if in 
every 
quantum game and for every quantum strategy $\p$, we 
have 
$P_1(\p,\mu)=P_1(\p,\mu')$ and 
$P_2(\p,\mu)=P_2(\p,\mu')$.

{\bf Example 1.4.}  The strategy supported on the 
singleton 
$\{\p\}$ is equivalent to the strategy supported on the 
singleton $\{-\p\}$ and to no other singleton.

{\bf Definition 1.5.}  Let $\nu$ be a mixed strategy and 
$\u$ 
a unit quaternion.  The {\sl right translate\/} of $\nu$ 
by 
$\u$ is the measure $\nu\u$ definied by $(\nu\u)(A)=\nu( 
A\u)$ where $A$ is any subset of the unit quaternions 
and 
$A\u =\{\x\u|x\in A\}$.  Similarly, the left translate 
of 
$\nu$ by \u is defined by $(\u\nu)(A)=\nu(\u A)$.  The 
following proposition is immediate:

{\bf Proposition 1.6.}  Let $(\nu,\mu)$ 
be a pair of mixed strategies and $\u$ a unit 
quaternion.  
Then in any game $\G^Q$, $(\nu,\mu)$ is a mixed strategy 
Nash 
equilibrium if and only if $(\nu\u,\u^{-1}\mu)$ is.

{\bf Definition 1.7.}  Two pairs of mixed strategies 
$(\nu,\mu)$ and $(\nu',\mu')$ are {\sl equivalent\/} if 
there 
exists a unit quaternion $\u$ such that $\nu'$ is 
equivalent 
to $\nu\u$ and $\mu'$ is equivalent to $\u^{-1}\mu$.  
Note 
that this definition is independent of any particular 
game.    

{\bf Proposition 1.8.}  In a given game, a pair of 
mixed 
strategies is a Nash equilibrium if and only if every 
equivalent pair of mixed strategies is also a Nash 
equilibrium.

\medskip

{\bf 2. Preliminary Results.}

Theorems 2.1, 2.2 and 2.4 are the main results which 
will be 
used in Section 3 to classify Nash equilibria.

{\bf Theorem 2.1.} 
Every mixed
strategy is equivalent to a mixed
strategy supported on (at most) four
points.  Those four points can be
taken to form an orthonormal basis for
$\R^4$.

{\bf Proof.}  First, choose any
orthonormal basis
$\{\q_1,\q_2,\q_3,\q_4\}$ for $\R^4$.
For any quaternion $\p$, write
(uniquely)$$\p=\sum_{\alpha=1}^4A_
\alpha(p)
\q_\alpha$$
where the $A_\alpha(p)$ are real
numbers.

Define a probability measure $\nu$
supported on the four points
$\q_\alpha$ by
$$\nu(\q_\alpha)=
\int_{\S^3}A_\alpha(\q)^2d\mu(\q)$$

For any two quaternions $\p$ and $\q$,
define
$$X(\p,\q)=\sum_{\alpha=1}^4
\pi_\alpha(\p)\pi_\alpha(\q)X_i
\eqno(2.1.1)$$

Then for any $\p$ we have
$$\eqalign{
P(\p,\mu)&=
\int_{\S^3}P(\p\q)d\mu(\q)\cr
&=
\int_{S^3}P\left(\sum_{\alpha=1}^4
A_\alpha(\q)\p\q_\alpha\right)
 d\mu(\q)\cr
&=
\sum_{\alpha=1}^4 P(\p\q_\alpha)
\int_{S^3}A_\alpha(\q)^2d\mu(\q)+
2\sum_{\alpha\neq\beta}
X(\p\q_\alpha,\p\q_\beta)
\int_{\S^3}A_\alpha(\q) A_\beta(\q)
d\mu(\q)\cr
&=P(\p,\nu)+ 2\sum_{\alpha\neq\beta}
X(\p\q_\alpha,\p\q_\beta)\int_{\S^3}A_
\alpha(\q) A_\beta(\q)
d\mu(\q)\cr}
$$

To conclude that $\mu$ is equivalent
to $\nu$ it is sufficient (and
necessary) to choose the $\q_\alpha$
so that for each $\alpha\neq\beta$ we
have
$$\int_{\S^3}A_\alpha(\q) A_\beta(\q)
d\mu(\q)=0$$

For this, consider the function
$B:\R^4\times\R^4\rightarrow\R$
defined by
$$B(\a,\b)=\int_{\S^3}
\pi_1(\overline{\a}\q)
\pi_1(\overline{\b}\q)d\mu(\q)$$
$B$ is a bilinear symmetric form and
so can be diagonalized; take the
$\q_\alpha$ to be an orthonormal basis
with respect to which $B$ is diagonal.
Then we have (for $\alpha\neq\beta$)
$$\eqalign{
\int_{\S^3}A_\alpha(\q)A_\beta(\q)
d\mu(\q)&=
\int_{\S^3}
\pi_1(\overline{\q_\alpha}\q)
\pi_1(\overline{\q_\beta}\q)
d\mu(\q)\cr
&=B(\q_\alpha,\q_\beta)=0\cr}$$
\qed

{\bf Theorem 2.2.}  Taking Player
2's (mixed) strategy $\mu$ as given,
Player 1's optimal response set is
equal to the intersection of $\S^3$
with a linear subspace of $\R^4$.

(Recall that we identify the unit
quaternions with the three-sphere
$\S^3$.)

{\bf Proof.}  Player One's problem is
to choose $\p\in S^3$ to maximize
$$P_1(\p,\mu)=\int P_1({\bf
pq})d\mu(\q)\eqno(2.2.1)$$

Expression (2.2.1) is a (real)
quadratic form in the coefficients
$\pi_i(\p)$ and hence is maximized
(over $S^3$) on the intersection of
$S^3$ with the real linear subspace of
$\R^4$ corresponding to the maximum
eigenvalue of that form.

\qed

{\bf Definition 2.3.}  We define the
function $K:\S^3\rightarrow\R$
by
$K(A+Bi+Cj+Dk)=ABCD$.
Thus in particular $K(\p)=0$ if and
only if \p is  a linear combination of
at most three of the fundamental units
$\{1,i,j,k\}$.

{\bf Theorem 2.4.}  Let $\mu$ be a mixed strategy 
supported 
on four orthogonal points $\q_1,\q_2,\q_3,\q_4$ played 
with 
probabilities $\alpha,\beta,\gamma,\delta$.  Suppose 
$\p$ is 
an optimal response to $\mu$ in some game where it is 
not the 
case that $X_1=X_2=X_3=X_4$.  Then $\p$ must satisfy:

$$\eqalign{&(\alpha-\beta)(\alpha-
\gamma)
    (\alpha-\delta)K(\p\q_1)+  (\beta-\alpha)(\beta-
\delta)(\beta-\gamma)K(\p\q_2)\cr
+& (\gamma-\alpha)(\gamma-\beta)
(\gamma-\delta)K(\p\q_3
)
+(\delta-\alpha)(\delta-\beta)
(\delta-
\gamma)K(\p\q_4)=0\cr}\eqno(2.4.1)
$$

{\bf Proof.}
Set $\p_n=\pi_n(\p)$ and consider the function
$$\matrix{
{\cal P}:&\S^3\times\R^4&\rightarrow&
\R\cr
&(\p,\x)&\mapsto &\sum_{n=1}^4
\p_n^2\x_n d\mu(\q)
\cr}$$ 

In particular, if we let
$X=(X_1,X_2,X_3,X_4)$
then ${\cal P}(\p,X)=P_1(\p,\mu)$.

The function ${\cal P}$ is quadratic
in \p and linear in \x; explicitly we
can write
$${\cal P}(\p,\x)=\sum_{i,j,k}
t_{ijk}\p_i\p_j\x_k$$ for some real
numbers $t_{ijk}$.

Set
$$M_{ij}(\x)=\sum_{k=1}^4 t_{ijk}
\x_k$$
$$N_{ij}(\p)=\sum_{k=1}^4 t_{ikj}
\p_j$$so that
$$M(\x)\cdot\pmatrix{\p_1\cr
\p_2\cr\p_3\cr\p_4\cr}
=
N(\p)\cdot\pmatrix{\x_1\cr
\x_2\cr\x_3\cr\x_4\cr}\eqno(2.4.2)$$

If $\p$ is an optimal response to the
strategy $\mu$, then
$(\p_1,\p_2,\p_3,\p_4)^T$ must be an
eigenvector of $M(X)$, say with
associated eigenvalue $\lambda$.
From this and (2.4.2) we conclude that
$$N(\p)\cdot\pmatrix{X_1\cr X_2 \cr
X_3 \cr X_4\cr}=
\lambda\cdot\pmatrix{
\p_1\cr \p_2\cr \p_3\cr\p_4\cr}
= N(\p)\cdot\pmatrix{\lambda\cr
\lambda\cr\lambda\cr\lambda\cr}$$where
the second equality holds by an
easy calculation.

Thus $N(\p)$ must be singular.  But it
follows from a somewhat less easy
calculation that the determinant of
$N(\p)/2$ is given by the left side
of (2.4.1).

\bigskip

{\bf 3.  Classification.}

{\bf Definition 3.1.}  Let \G be a two-by-two game with 
payoff pairs $(X_1,Y_1),\ldots,(X_4,Y_4)$ (listed in 
arbitrary order).  
$\G$ is a {\sl generic game\/} if the $X_i$ are all 
distinct, 
the $Y_i$ are all distinct, the twofold sums $X_i+X_j$ 
are 
all distinct and the twofold sums $Y_i+Y_j$ are all 
distinct.  

Theorem 3.3 will classify Nash Equilibria in $\G^Q$ 
where $\G$ is any generic two-by-two game.  Subtler 
versions 
of 
the same arguments work for non-generic games (yielding
somewhat messier results); see [NE].

To state Theorem 3.3 we need a definition:

{\bf Definition 3.2.}  Let $\p,\q,\r,\s$ be quaternions; 
write
$\p=p_1+p_2i+p_3j+p_4k$, etc.  Then the quadruple 
$(\p,\q,\r,\s)$ is {\sl intertwined} if there is a 
nonzero 
constant $\alpha$ such that 
$$\alpha(X\p+Y\q)=X\r+Y\s$$
identically in the polynomial variables $X$ and $Y$.  

Thus if the components of $\p,\q,\r,\s$ are all nonzero, 
then 
$(\p,\q,\r,\s)$ is intertwined if and only if the four 
quotients 
${p_1\over q_1},{p_2\over 
q_2},{p_3\over q_3},{p_4\over q_4}$ are equal (in some 
order) 
 to the four quotients ${r_1\over 
s_1},{r_2\over s_2},{r_3\over s_3},{r_4\over s_4}$.

The intertwined quadruple $(\p,\q,\r,\s)$ is {\sl fully 
intertwined} if $(\p,\r,\q,\s)$ is also intertwined.

We can now state the main theorem:

{\bf Theorem 3.3.} Let $\G$ be a generic game.  Then up 
to 
equivalence, every equilibrium in $\G^Q$ is of one of 
the 
following types:
\itemitem{a)}Each player plays each of four orthogonal 
quaternions with probability 1/4.
\itemitem{b)}Each player's strategy is supported on 
three of 
the four quaternions $1,i,j,k$.
\itemitem{c)}$\mu$ is supported on two orthogonal points 
$1,\v$; $\nu$ is supported on two orthogonal points 
$\p,\pu$, 
and the quadruple $(\p,\pv,\pu,\pvu)$ is fully 
intertwined.
\itemitem{d)}Each of $\mu$ and $\nu$ is supported on two 
orthogonal points, each played with probability $1/2$.  
Moreover, the supports of $\mu$ and $\nu$ lie in 
parallel 
planes.
\itemitem{e)}Each player plays a pure strategy from the 
four 
point set $\{1,i,j,k\}$.

{\bf Proof.}  Let $(\nu,\mu)$ be an equilibrium.  By 
(2.1) we 
can assume that each of $\nu$ and $\mu$ is supported on 
a set 
of at most four orthogonal points.  Applying a 
translation as 
in (1.7) we can assume that the support of $\mu$ 
contains 
the 
quaternion 1.  Then from standard facts about 
orthogonality 
in the space of quaternions, the support of $\mu$ is 
contained in a set of the form $\{1,\u,\v,\uv\}$ where 
$\u^2=\v^2=-1$ and $\uv+\vu=0$, played with 
probabilities of 
$\alpha,\beta,\gamma,\delta\ge0$.  We will maintain 
these 
assumptions and this notation while proving Theorems 
3.4, 
3.5, 3.9, and 3.10, which together imply Theorem 3.3.

{\bf Theorem 3.4.}  $\nu$ is a pure strategy if and only 
if 
$\mu$ is a pure strategy.  

{\bf Proof.}  If $\nu$ is a pure strategy, Player Two 
can 
guarantee any desired probability distribution over four 
outcomes; by genericity his optimal probability 
distribution 
is unique.\qed

{\bf Theorem 3.5.}  If the support of $\nu$ contains 
four 
points 
then $\mu$ assigns probability $1/4$ to each of four 
strategies.  

{\bf Proof.}  Explicitly write 
$\u=Ai+Bj+Ck,\v=Di+Ej+Fk,\u\v=Gi+Hj+Ik$.  
Write
$${\cal M}=\pmatrix{AB&DE&GH\cr AD&DF&GI\cr BC&EF&HI\cr 
}$$

By (2.2) the quadratic form 
$$\p\mapsto P_1(\p,\mu)\eqno(3.5.1)$$ is constant on the 
unit 
sphere $\S^3$.  Therefore its non-diagonal coefficients 
are 
all zero.  Computing these coefficients explicitly and 
dividing by (non-zero) expressions of the form $(x_i-
x_j)$, 
we get 
$${\cal 
M}\cdot(\beta,\gamma,\delta)^T=(0,0,0)^T\eqno(3.5.2)$$
But ${\cal M}$ also kills the column vector $(1,1,1)^T$.  
Thus we have two cases:

Case I.  $\beta=\gamma=\delta$.  Then the four diagonal 
terms 
of (3.5.1) (which must all be equal) are given   by 
$(X_1+X_2+X_3+X_4)\beta+X_i(\alpha-\beta)$, with 
$i=1,2,3,4$.  
Since the $X_i$ are not all equal, it follows 
$\alpha=\beta=\gamma=\delta=1/4$, proving the theorem.

Case II.  {\cal M} has rank at most one.  From this and 
the 
orthogonality of \u,\v,\u\v, we have 
$\{\u,\v,\u\v\}\cap\{i,j,k\}\neq\emptyset$.  Assume 
$\u=i$ 
(the other cases are similar).  Then $A=1$, $B=C=D=G=0$, 
$H=-
F$ and $I=E$.   
The four diagonal entries of (3.3.1) are now equal; call 
their common value $\lambda$ so that we have
$$\pmatrix{
\alpha&\beta&E^2\gamma+F^2\delta&E^2\delta+F^2\gamma\cr
\beta&\alpha&E^2\delta+F^2\gamma&E^2\gamma+F^2\delta\cr
E^2\gamma+F^2\delta&E^2\delta+F^2\gamma&\alpha&\beta\cr
E^2\delta+F^2\gamma&E^2\gamma+F^2\delta&\beta&\alpha\cr}
\cdot
\pmatrix{X_1\cr X_2\cr  X_3\cr  X_4\cr}=
\pmatrix{\lambda\cr\lambda\cr\lambda\cr\lambda\cr}\eqno(
3.5.3
)$$
Combining (3.5.2), (3.5.3), the conditions 
$\alpha+\beta+\gamma+\delta=E^2+F^2=1$ and the 
genericity 
conditions, we get $\alpha=\beta=\gamma=\delta$ as 
required.\qed

{\bf Corollary 3.5A.} If either player's strategy
has a four-point support, then each player plays each of
four
orthogonal quaternions
with probability $1/4$.

{\bf Proof.}  Apply Theorem 3.5 twice, one
as
stated and once with the players reversed.\qed

Theorem 3.9, dealing with the case where $\nu$ is 
supported 
on exactly three points, requires some preliminary 
lemmas:

{\bf Lemma 3.6.}  It is not the case that Player Two 
plays 
$1,\u,\v$ each with probability $1/3$.

{\bf Proof.}  If $1,\u,\v$ are played with probability 
$1/3$
then one computes that the eigenvalues of the form 
(2.2.1) 
are $X_1+X_2+X_3$, $X_1+X_2+X_4$, $X_1+X_3+X_4$, 
$X_2+X_3+X_4$, which are all distinct by genericity.  
Thus 
Player One responds with a pure strategy, and Theorem 
3.4 
provides a contradiction.
\qed

{\bf Lemma 3.7.}  Suppose the support of $\nu$ is 
contained 
in the linear span of $1,i,j$. and suppose that $1$ and 
$i$ 
are both optimal responses for Player Two.  Then one of 
the 
following is true:
\itemitem{a)}The support of $\nu$ is contained in the 
three 
point set $\{1,i,j\}$
\itemitem{b)}The support of $\nu$ is contained in a set 
of 
the form $\{1,Ei+Fj,-Fi+Ej\}$ with $Ei+Fj$ and $-Fi+Ej$ 
played equiprobably.

Moreover, if b) holds and either $j$ or $k$ is also an 
optimal response for Player Two, then 1 is played with 
probability zero.

{\bf Proof.}  Suppose $\nu$ is 
supported on three orthogonal quaternions $\q_1=
A+Bi+Cj$, $\q_2=D+Ei+Fj$, $\q_3=G+Hi+Ij$, played with 
probabilities $\phi,\psi,\xi$.  The first order 
conditions 
for Player Two's maximization problem must be satisfied 
at 
both $1$ and $i$; this (together with genericity for the 
game 
$\G$) gives
$$\pmatrix{AC&DF&GI\cr BC&EF&HI\cr}\pmatrix{
\phi\cr \psi\cr\xi\cr}=
\pmatrix{0\cr0\cr0\cr}=
\pmatrix{AC&DF&GI\cr BC&EF&HI\cr}\pmatrix{
1\cr1\cr1\cr}\eqno(3.7.1)$$
so that 
by (3.6) with the players reversed, the matrix on the 
left 
has rank at most 
one.  
This (together with the orthogonality of 
$\q_1,\q_2,\q_3)$ 
gives $\{\q_1,\q_2,\q_3\}\cap\{1,i,j\}\neq
\emptyset$.  We can assume $\q_1=1$ (all other cases are 
similar); thus $A=1$, $B=C=D=G=0$,  $H=-F$, $I=E$.  Now 
(3.7.1)
says $EF(\psi-\xi)=0$.  If $EF=0$, then a) holds and if 
$\psi-\xi=0$ then b) holds.  

Now suppose $j$ is also an optimal response for Player 
Two. 
Then $0=P_2(\nu,i)-P_2(\nu,j)=\phi(Y_2-
Y_3)$, so that 
by genericity $\phi=0$.  A similar argument works if $k$ 
is 
optimal.\qed

{\bf Lemma 3.8.} Suppose $\nu$ is supported on exactly 
three 
points and continue to assume that $\mu$ is supported on 
a 
subset of $\{1,\u,\v\,\uv\}$.  Then at least two of the 
four 
quaternions $1,\u,\v,\uv$ are optimal responses for 
Player 
One.

{\bf Proof.}  By (3.5), $\mu$ is supported on at most 
three 
points; we can rename so those points are $1,\u,\v$.  
These 
are played with probabilities $\alpha,\beta,\gamma$ and 
we 
can rename again so that $\alpha$ lies (perhaps not 
strictly) 
between $\beta$ and $\gamma$.

If $\p$ is any optimal response by Player One, apply 
(2.4) 
with $\delta=0$ (and possibly $\gamma=0$) to get 
$$\sigma_1K(\p)+\sigma_2K(\p\u)+\sigma_3K(\p\v)+\sigma_4
K(\p
\u\v)=0\eqno(3.8.1)$$
where $\sigma_1=(\alpha-\beta)(\alpha-\gamma)\alpha$, 
etc., 
so that  
$$\sigma_1,\sigma_4\le0\qquad\hbox{and}\qquad\sigma_2,
\sigma_3
\ge0\eqno(3.8.2)$$

Case I:  Suppose none of the  $\sigma_i$ is equal to 
zero.  Then $\gamma\neq0$ so a) holds.

By (2.2), the support of $\nu$ spans a three-dimensional 
hyperplane in $\R^4$ and thus must include some 
quaternion of 
the form $A+B\u$ ($A,B\in\R$).   Inserting $\p=A+B\u$ 
into 
(3.8.1) gives
$$AB(\sigma_1 B^2-\sigma_2 A^2)K(1+\u)=0\eqno(3.8.3)$$ 
Thus either $AB=0$ (in which case either $\p=1$ 
or $\p=\u$) or $K(1+\u)=0$.  This and similar arguments 
establish
the following:
$$\hbox{If $1$ and $\u$ are both suboptimal responses, 
then 
$K(1+\u)=0$.}\eqno(3.8.3a)$$
$$\hbox{If $1$ and $\v$ are both suboptimal responses, 
then 
$K(1+\v)=0$.}\eqno(3.8.3b)$$
$$\hbox{If $\u$ and \uv are both suboptimal responses, 
then 
$K(1+\v)=0$.}\eqno(3.8.3c)$$
$$\hbox{If $\v$ and $\uv$ are both suboptimal responses, 
then 
$K(1+\u)=0$.}\eqno(3.8.3d)$$

Taken together, these imply that if the lemma fails, 
then 
$K(1+\u)=K(1+\v)=0$.  From this it follows that 
$\{\u,\v,\uv\}\cap\{\pm i,\pm j, \pm k\}\neq\emptyset$; 
assume without loss of 
generality that $\u=i$ and therefore $\v$ is in the 
linear 
span of $\{j,k\}$.  (Generality is not lost because the 
argument to follow works just as well, with obvious 
modifications, in all the remaining cases.)

Now we have 
$$\eqalign{
P_1(A+B\v,\mu)&=\alpha P_1(A+B\v)+ 
\beta P_1(Ai+B\v i)+\gamma P_1(A\v-B)\cr
&=
A^2\Big(\alpha P_1(1)+\beta P_1(i)+\gamma P_1(\v)\Big) +
B^2\Big(\alpha P_1(\v)+\beta P_1(\v i)+\gamma 
P_1(1)\Big)}$$
which is maximized at an endpoint, so either 1 or $\v$
is an optimal response for Player One.  
Similarly, at least one from each pair $\{1,\uv\}$, 
$\{\u,\v\}$, and $\{\u,\uv\}$ is an optimal response, 
from 
which b) (and therefore the lemma) follows.

Case II:  Suppose at least one of the $\sigma_i$ is 
equal to 
zero.  Up to renaming $\u$ and $\v$, there are three 
ways 
this can happen:

Subcase IIA:  $\alpha=\beta$, $\gamma=0$.  As above, 
Player 
One's optimal response set contains a quaternion of the 
form 
$(A+B\u)$.  But $P_1(A+B\u,\mu)$ is independent of $A$ 
and 
$B$, so both $1$ and $\u$ are optimal, proving the 
theorem.  
(Note that $\v$ and $\uv$ are also both optimal, so that 
in 
fact by (3.5A) this case never occurs.) 

Subcase IIB:  $\alpha=\beta$, $\gamma\neq0$.  By Lemma 
(3.6), 
$\gamma\neq\alpha,\beta$.  Thus $\sigma_3$ and 
$\sigma_4$ are 
nonzero, so (3.8.3b), (3.8.3c) and (3.8.3d) (but not 
(3.8.3a)) still hold.  But $\sigma_1=\sigma_2=0$ so the 
same 
techniques now yield
$$\hbox{If $1$ and $\u$ are both suboptimal responses, 
then 
$K(1+\u)=0$.}\eqno(3.8.3e)$$
$$\hbox{If $1$ and $\v$ are both suboptimal responses, 
then 
$K(1+\v)=0$.}\eqno(3.8.3f)$$
We can now repeat the argument from Case I.

Subcase IIC:  $\alpha\neq\beta$, $\gamma=0$.  Now we 
have 
$\sigma_1,\sigma_2\neq0$, $\sigma_3=\sigma_4=0$, so that 
(3.8.3a) through (3.8.3c) still hold, along with 
(3.8.3e) and 
(3.8.3f).  We can now repeat the argument from Case I.
\qed

{\bf Theorem 3.9.}  If $\nu$ is supported on exactly 
three 
points, then up to equivalence, both $\mu$ and $\nu$ are 
supported on three-point subsets of $\{1,i,j,k\}$.

{\bf Proof.}  By (3.5) we can assume that $\mu$ is 
supported 
on $\{1,\u,\v\}$.  By (3.8) we can assume without much 
loss 
of generality that $1$ and $\u$ are optimal responses 
for 
Player One.  (The argument below works equally well, 
with 
obvious modifications, for other pairs.)  Let $\w$ be a 
quaternion orthogonal to 1 and \u such that the support 
of 
$\nu$ is contained in the linear span of $1,\u$ and 
$\w$. 

By (2.2), any quaternion of the form $X+Y\u+Z\w$ is an 
optimal response for Player One, so by (2.4) we have
$$
\sigma_1 K(X+Y\u+Z\w)+
\sigma_2 K(X\u-Y+Z\w\u)+
\sigma_3 K(X\v+Y\u\v+Z\w\v)+
\sigma_4 K(X\u\v-Y\v+Z\w\u\v)=0$$
identically in $X,Y,Z$.  Writing out the left side as a 
polynomial in these three variables, the coefficients, 
all of 
which must vanish, can be expressed in terms of the 
components of $\u,\v,\w$.  Setting all these expressions 
equal to zero and solving, we find that 
$\{\u,\v,\w\}\in\{\pm 
i,\pm j,\pm k\}$.  (The details of this tedious but 
straightforward calculation can be found on pages 32-33 
of 
[NE].)   We assume $\u=i$, $\w=j$.  

Claim:  Player Two's strategy is not supported just on 
$1$ 
and $i$.  Proof:  If so, the fact that 
$P_1(1,\mu)=P_1(i,\mu)$ implies that $\mu$ assigns equal 
weights to $1$ and $i$, which implies 
$P_1(j,\mu)=P_1(k,\mu)$, contradicting the fact that $j$ 
but 
not $k$ is optimal for Player One.

Thus the support of $\mu$ is a three-point subset of 
$\{1,i,j,k\}$.  It now follows from Lemma (3.8) 
(together 
with the assumption that the support of $\nu$ contains 
three 
points) that the support of $\nu$ is $\{1,i,j\}$, 
completing 
the proof.\qed

{\bf Theorem 3.10.}  Suppose $\nu$ is supported on two 
points.  
Then $\mu$ is supported on $1,\u$ and $\nu$ is supported 
on 
two quaternions $\p,\pv$  where either
\itemitem{a)}The quadruple $(\p,\pu,\pv,\pvu)$ is fully 
intertwined or
\itemitem{b)}$\u=\v$ and each player plays each strategy 
with 
probability $1/2$.

{\bf Proof.}  Suppose $1$ and $\u$ are played with 
probabilities $\alpha$ and $\beta$.  

Any unit quaternion of the form $X\p+Y\pv$ is an optimal 
response for Player One; thus (2.4) with 
$\q_1=1,\q_2=\u,\gamma=\delta=0$ gives
$$(\alpha-\beta)\Big(\alpha^2K(X\p+Y\pv)-
\beta^2K(X\pv+Y\pvu)\Big)=0$$.  This, plus the identical 
observation with the players reversed, estabilishes full 
intertwining except when $\alpha=\beta=1/2$.  In that 
case, 
$P_1(\p,\mu)=P_1(\pu,\mu)$ so $\pu$ must be optimal; 
i.e. we 
can take $\v=\u$.\qed

This completes the proof of Theorem 3.3.

{\bf Remark.}  The statement of Theorem 3.3 makes it 
natural to ask for a classification of fully intertwined 
quadruples of the form $(\p,\pv,\pu,\pvu)$ with $\u$, 
$\v$ square roots of $-1$.  That classification is 
provided in [I].  The thrust of the result is this:  All 
such quadruples fall into one of approximately 15 
families.  Each of these families is at most 
four-dimensional (inside the twelve-dimensional
manifold of all four-tuples).  For all but one of the families, it is 
easy to tell by inspection whether a given quadruple 
satisifies the membership condition.  The exceptional 
family is one-dimensional.

In short:  Condition b) of Theorem 3.3 allows only four 
dimensions worth of possible equilibria, all of which 
are easily identifiable except for a one-dimensional 
subset.

\bigskip

{\bf 4.  Minimal Payoffs and Opting Out}

Theorem 3.3. classifies all mixed strategy Nash 
equilibria in 
generic games.   Here we briefly address the issue of 
whether 
these equilibria survive in a larger game where the 
players 
can opt out of  the assigned communication protocol.  

A key tool is the very simple Theorem 4.1; this and its 
corollary 4.1A apply to all two by two games (whether 
generic 
or not) and are of independent interest:

{\bf Theorem 4.1.}  Let $\G$ be a game with payoff pairs 
$(X_1,Y_1),\ldots(X_4,Y_4)$.  Then in any mixed strategy 
quantum equilibirum, Player One earns at least 
$(X_1+X_2+X_3+X_4)/4$.  

{\bf Proof.}  Player One maximizes the quadratic form 
(2.2.1) 
over the sphere $\S^3$.  The trace of this form is 
$X_1+X_2+X_3+X_4$, so the maximium eigenvalue must be at 
least $(X_1+X_2+X_3+X_4)/4$.\qed

{\bf Corollary 4.1A.}  If, in Theorem 4.1, the game $\G$ 
is 
zero-sum, then in any mixed strategy quantum 
equilibrium, 
Player One earns exactly $(X_1+X_2+X_3+X_4)/4$.

{\bf Proof.}  Apply (4.1) to both players.\qed

{\bf 4.2. Remarks on Opting Out.}  A player can 
throw away his entangled penny and substitute an 
unentangled 
penny (or for that matter a purely classical penny, but 
this 
offers no additional advantage, because the unentangled 
quantum penny can always be returned in one of the two 
classical states $\H$ or $\T$).  However, a simple 
quantum 
mechanical calculation shows that if Player One 
unilaterally 
substitutes an unentangled penny, then no matter what 
strategies the players follow from there, the result is 
a 
uniform distribution over the four possible outcomes.  
By 
Theorem 4.1, Player One considers this weakly inferior 
to 
any $\G^Q$ equilibrium.  Thus, even if we allow players 
to 
choose their pennies, all of the $\G^Q$ equilibria 
survive.

\bigskip

\centerline{\bf References}

\item{[A]}R.~J.~Aumann, ``Subjectivity and Correlation 
in 
Randomized Strategies'', {\it J. of Mathematical 
Economics} 1 
(1974).

\item{[CHTW]}R.~Cleve, P.~Hoyer, B.~Toner and 
J.~Watrous, 
``Consequences and Limits of Nonlocal Strategies'', {\it 
Proc. of the 19th Annual Conference on Computational 
Complexity\/} (2004), 236-249.

\item{[DL]}G.~Dahl and S.~Landsburg, ``Quantum 
Strategies in 
Noncooperative Games'', preprint, available at
\centerline{http://www.landsburg.com/dahlcurrent.pdf}

\item{[EW]}J.~Eisert and M.~Wilkens, ``Quantum Games'', 
{\it 
J. of Modern Optics\/} 47 (2000), 2543-2556

\item{[EWL]}J.~Eisert, M.~Wilkens and M.~Lewenstein,  
``Quantum Games and Quantum Strategies'', {\it Phys. 
Reve. 
Lett. 83} (1999), 3077-3080

\item{[I]}S.~Landsburg, ``Intertwining'', working paper 
available at 

\centerline{http://www.landsburg.com/twining4.pdf}

\item{[M]}P.~La Mura, ``Correlated Equilibria of 
Classical 
Strategic Games with Quantum Signals'':
\centerline{arXiv:quant-ph/0390933v1}, September 2003

\item{[NE]}S.~Landsburg, ``Nash Equilibrium in Quantum 
Games'', RCER Working Paper \#524, Rochester Center for 
Economic Research, February 2006

\bye